# FUZZY MATHEMATICAL MODEL FOR OPTIMIZING SUCCESS CRITERIA OF PROJECTS: A PROJECT MANAGEMENT APPLICATION

*Mohammad Sammany[1], Ahmad Steef [2*], Nedaa Agami[3], T. Medhat[4]*

[1] *Faculty of Sciences, Department of Mathematics, University of Aleppo, Syria,*
[2*] *Faculty of Informatics Engineering, University of Idlib, Syria,*
,[3] *Faculty of Computers and Artificial Intelligence, Operations Research Department, Cairo University, Cairo*
[4] *Department of Electrical Engineering, Faculty of Engineering, Kafrelsheikh University, Kafrelsheikh 33516, Egypt*

*Dr.sammany@hotmail.com, msma7766@gmail.com, nedaa.agami@gmail.com, tmedhatm@eng.kfs.edu.eg*

*[*]Corresponding Author: msma7766@gmail.com*

*Abstract*— It is well known over the recent years that measuring the success of projects under the umbrella of project management is inextricably linked with the associated cost, time, and quality. Most of the previous researches in the field assigned a separate mathematical model for each criterion, then numerical methods or search techniques were applied to obtain the optimal trade-off between the three criteria. However in this paper, the problem was addressed by linear multi-objective optimization using only one fuzzy mathematical model. The three criteria were merged in a single non-linear membership function to find the optimal trade-off. Finally, the proposed model is tested and validated using numerical examples.

## I. INTRODUCTION AND

To ensure the success of a project from the beginning to end, it is crucial to take into consideration that all expenditures are covered by available finances during the execution time. However, the cost of a project is strongly correlated to the entire execution time, which can be controlled either by increasing the number of resources (labor) or working hours, which means paying additional costs. Therefore, minimizing the project time and cost is the main objective for the success of projects from a project management perspective. In addition, the instant determination of expenditures and comparing them to the project budget is very beneficial in controlling entire costs of a project according to the specified time period. On the other hand, a quality criterion which describes the degree of execution efficiency is also important for the success of a project. Due to its inextricable link with time and cost, the quality criterion can be very critical to the project if there are no adjustable restrictions on it. Most of the previous researches tackling this problem assigned one mathematical model for each criterion, then numerical methods or search techniques were applied to obtain the optimal trade-off between the three criteria. For example, in [1], Babo et.al presented three mathematical models for the associated cost, time and quality of a project by assuming that the cost and quality of an activity vary as linear functions with the execution time. Another study presented by D. Khang et.al [2], which has a significant impact on industrial applications, is also based on [1]. Both results in [1] and [2] were applied on AOA networks for projects [3]. On the other hand, in [4] Elrayess and Kandil used a multi-objective genetic algorithm to optimize cost, time and quality of a project, whereas in [5], [6] a multi

Ant Colony (AC) algorithm is used to solve the same problem. Accordingly, we raise the following question: Can this problem be solved using only one mathematical model and give the same results of previous studies while helping the decision makers to adjust the cost and time of a project simultaneously according to proposed constraints related to any of the criteria? The aim of this research is to try to find an answer to this question. In other words, to identify connection between cost, time and quality criteria for a project management problem, besides providing the best scenarios for project planners. Further to that, allow for controlling the execution process in such a way that satisfies the required objectives within the available capabilities. For this purpose, three mathematical models describing the project time and cost and quality were proposed. Based on Critical Path Method (CPM), execution time of the project was calculated and analyzed. Given that the cost and quality of an activity vary as linear functions during the execution time, the proposed models were merged into a single fuzzy mathematical model using non-linear membership function to find the optimal trade-off between the three objectives as a solution of linear multi-objective optimization problem (MOOP) [7]. Finally the proposed model is demonstrated and tested using a numerical example, and the optimal solutions are obtained.

The rest of the paper is organized as follows: Section 2 consists of 3 subsections where the first one discusses the background needed for constructing the model, while the second and the third subsections present the membership function and the fuzzy mathematical model respectively. Section 3 is dedicated to the empirical study and illustration of the proposed method using a numerical example. Finally, the paper ends with the conclusion and suggestions for future work.

## II. METHODOLOGY

Suppose that there is a linear relationship between cost, time, and quality for a project with $n$ activities. Let $c_j$ be the normal cost of executing activity $j$, which corresponds to the normal time $t_j$ and the normal quality $q_j$, and $q'_j$ is the crash quality that corresponds to crash cost $c'_j$ and crash time $t'_j$. Let $Q$, $T$, and $C$ be the straight lines of the linear functions for quality, time, and cost, respectively (see fig. 1). Then we can write:

$$C = c_j + \frac{c'_j - c_j}{t'_j - t_j}(T - t_j) \quad (1)$$

$$Q = q_j - \frac{q_j - q'_j}{t_j - t'_j}(t_j - T) \quad (2)$$

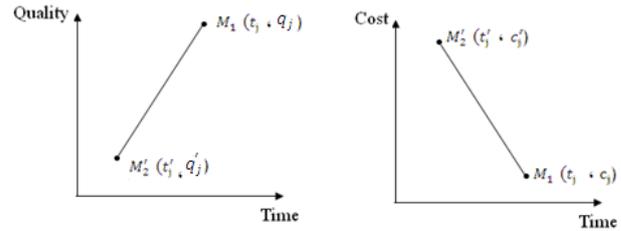

**Figure (1)**: Linear functions for the changes in cost, time, and quality of an activity

Before using these equations in the model, we give a glance on membership function in the next subsection.

### 2.2 Membership Function

Let $Z_k$; $k = 1, 2, ..., n$ be objective functions of the MOOP, and $X^k$ be their initial solution. This can be represented in a matrix of the form:

$$\begin{array}{c} \\ X^{(1)} \\ X^{(2)} \\ \vdots \\ X^{(k)} \end{array} \begin{bmatrix} Z_1 & Z_2 & \cdots & Z_k \\ Z_{11} & Z_{12} & \cdots & Z_{1k} \\ Z_{21} & Z_{22} & \cdots & Z_{2k} \\ \vdots & \vdots & \vdots & \vdots \\ Z_{k1} & Z_{k2} & \cdots & Z_{kk} \end{bmatrix} \quad (4)$$

where:

$$\left. \begin{array}{l} Z_{ij} = Z_j(X^{(i)}) \\ U_k = (Z_k)_{\max} = \max\{Z_{1k}, Z_{2k}, ..., Z_{kk}\} \\ L_k = (Z_k)_{\min} = \min\{Z_{1k}, Z_{2k}, ..., Z_{kk}\} \end{array} \right\} \quad (5)$$

In general, the approaches for converting multi-objectives into a single objective optimization problem has been recently addressed by many researchers using different types of membership functions such as linear, exponential, hyperbolic,…,etc. (see for instance [8], [9], [10]).

In [11], the authors presented mathematical model by liner membership function to control cost and time of a project, and in [12], the authors also presented mathematical model by linear membership function to minimize the cost and time for the Solid Transportation Problem, while in [13], the authors presented mathematical model by exponential membership function for controlling the cost and time of a project.

Since the problem understudy is concerned with finding the solutions that minimize the entire cost and time, and maximize quality of a project simultaneously, the purpose of using membership function here is to describe any of the three criteria at sufficiently high membership degree [14], [15] apart from its capability in controlling the solutions of linear multi-objective optimization. In the research work demonstrated in [16], [17], the authors used hyperbolic membership function given by the equation,

$$\mu_{Z_k}(x) = \begin{cases} 1 & \text{if } Z_k \leq L_k \\ \frac{1}{2}\left(Th\left(\left(\frac{U_k+L_k}{2} - Z_k(x)\right)a_k\right)\right) + \frac{1}{2} & \text{if } U_k < Z_k < L_k \quad (6) \\ 0 & \text{if } Z_k \geq U_k \end{cases}$$

where $a_k(x) = \dfrac{6}{U_k - L_k}$ ;

$k = 1, 2, ..., K$ \hspace{2cm} (7)

to solve minimization problems. Accordingly, the same function was used in the problem understudy.

### III. RESULTS AND DISCUSSION

**Fuzzy Mathematical Model**
Since we are solving a project management problem with the objective of minimizing the project execution time, we built our model on AON network utilized from CPM [18], [19], as follows:

$$\begin{cases} Min \; Z_1 = Y_{End} \\ \underline{\underline{S.T}} \\ \quad t'_j \leq T_j \leq t_j \; ; \quad j = 1,...,n \\ \quad Y_i + T_i \leq Y_j \; ; \quad i, j = 1,...,n \end{cases} \quad (8)$$

$Y_{End}$ : Start time of the last node in a network (end node) that doesn't have time

$T_j$ : New times of activities calculated according to the new time of the project

$Y_j$ : Start time of activity $j$

$Y_i$ : Start time of every activity $i$ that directly precedes activity $j$

In a similar manner, the mathematical models for total cost and quality can be easily written respectively as follows:

$$\begin{cases} Min \; Z_2 = \sum_{j=1}^{n}\left(c_j + \left(\dfrac{c'_j - c_j}{t_j - t'_j}\right)(t_j - T_j)\right) \\ \underline{\underline{S.T}} \\ \quad t'_j \leq T_j \leq t_j \; ; \quad j = 1,...,n \\ \quad Y_i + T_i \leq Y_j \; ; \quad i, j = 1,...,n \end{cases} \quad (9)$$

$$\begin{cases} MaxZ3 = \sum_{j=1}^{n}\left(q_j - \left(\dfrac{q_j - q'_j}{t_j - t'_j}\right)(t_j - T_j)\right) \\ \underline{\underline{S.T}} \\ \quad t'_j \leq T_j \leq t_j \; ; \quad j = 1,...,n \\ \quad Y_i + T_i \leq Y_j \; ; \quad i, j = 1,...,n \end{cases} \quad (10)$$

This can also be re-written as:

$$\begin{cases} MinZ3 = \sum_{j=1}^{n}\left(\dfrac{q_j - q'_j}{t_j - t'_j}\right)(t_j - T_j) \\ \underline{\underline{S.T}} \\ \quad t'_j \leq T_j \leq t_j \; ; \quad j = 1,...,n \\ \quad Y_i + T_i \leq Y_j \; ; \quad i, j = 1,...,n \end{cases} \quad (11)$$

Injecting the previous models (8), (9) and (11) in the non-linear membership function given by equation (6), we get a fuzzy mathematical model in the form:

$$\begin{cases} Max \; \lambda \\ \underline{\underline{ST}} \\ \lambda \leq \mu_{Z_1}(X) \\ \lambda \leq \mu_{Z_2}(X) \\ \lambda \leq \mu_{Z3}(X) \\ 0 \leq \lambda \\ t'_j \leq T_j \leq t_j \quad ; j = 1,...,n \\ Y_i + T_i \leq Y_j \end{cases} \quad (12)$$

Although our approach can flexibly be applied to include many cases of time, cost, and quality, we confined our study on some of these cases. To illustrate the idea, we demonstrate a numerical example in the next section.

### IV. EMPIRICAL STUDY

**Numerical Examples**
Using the optimization tool provided by Linear Interactive and Discrete Optimization program (LINDO), a code was written to solve the problem utilized from the MAX sub-routine (see Appendix).

Suppose we have an executive project with 9 activities. The normal cost and normal time, and changes of cost, time, and quality in terms of the crash time and crash cost for each activity are listed in table (1). The problem is solved by finding the optimal time, cost, and quality of the project in the following steps:

**Table (1)**: Dependencies, cost, time, and quality for each activity

| Activity | Dependency | Normal Time | Normal Cost | Crash Time | Crash Cost | Crash Quality |
|---|---|---|---|---|---|---|
| A | NON | 10 | 500000 | 6 | 700000 | 85% |
| B | A | 9 | 450000 | 7 | 600000 | 88% |
| C | B | 4 | 150000 | 3 | 210000 | 95% |
| D | B | 6 | 120000 | 4 | 200000 | 80% |
| E | B | 3 | 300000 | 2 | 400000 | 82% |
| F | C , D | 7 | 210000 | 5 | 290000 | 93% |
| G | D , E | 5 | 400000 | 3 | 550000 | 90% |
| H | E | 6 | 330000 | 3 | 510000 | 75% |
| I | F , G , H | 10 | 600000 | 7 | 840000 | 80% |

- **Step 1.** Formulate the fuzzy mathematical model:

$$\left. \begin{array}{l} Max\ \lambda \\ \underline{\underline{ST}} \\ \lambda \leq \frac{1}{2}Th\left(\left(\frac{U_1+L_1}{2}-Z_1(x)\right)a_1\right)+\frac{1}{2} \\ \lambda \leq \frac{1}{2}Th\left(\left(\frac{U_2+L_2}{2}-Z_2(x)\right)a_2\right)+\frac{1}{2} \\ \lambda \leq \frac{1}{2}Th\left(\left(\frac{U_3+L_3}{2}-Z_3(x)\right)a_3\right)+\frac{1}{2} \\ t'_j \leq T_j \leq t_j \quad ; j=1,...,9 \\ Y_i + T_i \leq Y_j \quad j=1,...,9 \\ \lambda \geq 0 \end{array} \right\} \quad (13)$$

- **Step 2.** Solve the first model in the form:

$$\left. \begin{array}{l} \min Z_1 = \sum_{j=1}^{9}\left(c_j + \left(\frac{c'_j - c_j}{t_j - t'_j}\right)(t_j - T_j)\right) \\ \underline{\underline{S.T}} \\ t'_j \leq T_j \leq t_j \ ; \quad j=1,...,9 \\ Y_i + T_i \leq Y_j \ ; \ i,j=1,...,9 \end{array} \right\} \quad (14)$$

It is easy to know that the obtained solution will be $Z_1(X^{(1)}) = 3060000$ at:

$X^{(1)} = \{T_A=10, T_B=9, T_C=4, T_D=6, T_E=3, T_F=7, T_G=5, T_H=6, T_I=10\}$

- **Step 3.** Solve the second model in the form:

$$\left( \begin{array}{l} Min\ Z_2 = Y_{End} \\ \underline{\underline{S.T}} \\ t'_j \leq T_j \leq t_j \ ; \quad j=1,...,9 \\ Y_i + T_i \leq Y_j \ ; \ i,j=1,...,9 \end{array} \right) \quad (15)$$

And similarly, the obtained solution will be $Z_2(X^{(2)}) = 29$ at:

$X^{(2)} = \{T_A=6, T_B=7, T_C=3, T_D=4, T_E=2, T_F=5, T_G=3, T_H=3, T_I=7\}$

- **Step 4.** Solve the third model in the form:

$$\left( \begin{array}{l} MinZ3 = \sum_{j=1}^{9}\left(\frac{q_j - q'_j}{t_j - t'_j}\right)(t_j - T_j) \\ \underline{\underline{S.T}} \\ t'_j \leq T_j \leq t_j \ ; \quad j=1,...,9 \\ Y_i + T_i \leq Y_j \ ; \ i,j=1,...,9 \end{array} \right) \quad (16)$$

The obtained solution: $Z_3(X^{(3)}) = 0$ at:

$X^{(3)} = \{T_A=10, T_B=9, T_C=4, T_D=6, T_E=3, T_F=7, T_G=5, T_H=6, T_I=10\}$

- **Step 5.** Compute $U_1, L_1, U_2,$ and $L_2$ based on the following matrix, e.q (4):

|  | $Z_1$ | $Z_2$ | $Z_3$ |
|---|---|---|---|
| $X^{(1)}$ | $Z_1(X^{(1)})$ | $Z_2(X^{(1)})$ | $Z_3(X^{(1)})$ |
| $X^{(2)}$ | $Z_1(X^{(2)})$ | $Z_2(X^{(2)})$ | $Z_3(X^{(2)})$ |
| $X^{(3)}$ | $Z_1(X^{(3)})$ | $Z_2(X^{(3)})$ | $Z_3(X^{(3)})$ |

$L_1 = 3060000, U_1 = 4250000 \Rightarrow Z_1 \in [3060000, 4250000]$
$L_2 = 29, U_2 = 42 \Rightarrow Z_2 \in [29, 42]$
$L_3 = 0, U_3 = 1 \Rightarrow Z_3 \in [0,1]$

- **Step 6.** Substitute the values computed in step 5 in the fuzzy mathematical model shown in step1.

$$\left. \begin{array}{l} Max\ \lambda \\ \underline{\underline{ST}} \\ \lambda \leq \frac{1}{2}Th\left(\left(\frac{3060000+4250000}{2}-Z_1\right)\left(\frac{6}{4250000-3060000}\right)\right)+\frac{1}{2} \\ \lambda \leq \frac{1}{2}Th\left(\left(\frac{42+29}{2}-Z_2\right)\left(\frac{6}{42-29}\right)\right)+\frac{1}{2} \\ \lambda \leq \frac{1}{2}Th\left(\left(\frac{0+1}{2}-Z_3\right)\left(\frac{6}{1}\right)\right)+\frac{1}{2} \\ t'_j \leq T_j \leq t_j \quad ; j=1,...,9 \\ Y_i + T_i \leq Y_j \quad j=1,...,9 \\ \lambda \geq 0 \end{array} \right\} \quad (17)$$

**The optimal solution is:**
$\lambda = 0.7997312\ ;\ Z_1 = 3440000\ ,\ Z_2 = 34\ , Z_3 = 0.34$
at
$X = \{T_A=6, T_B=7, T_C=4, T_D=6, T_E=3, T_F=5, T_G=5, T_H=6, T_I=10\}$.

This solution gives a balanced initial optimal solution for the three criteria because the model was built to search for some solution like $X$ such that:
$Z_1(X) \in [3060000, 4250000], Z_2(X) \in [29,42], Z_3(X) \in [0,1]$, and do minimize $Z_1, Z_2$ at the same time while maximizing $Z_3$. In other words, the proposed model conducts many tests to find the optimal solution for three criteria at the same time.

After the initial optimal solution has been obtained, decision makers can control the optimal solution of cost, time, and quality (versus each other) by adding constraints related to each criterion to the model (17).
For example, suppose we want to control the quality of some critical activities, which are important for us, so we will formulate constraints describing the quality of activities. It is easy to see that the constraint:

$$q_j - \frac{q_j - q'_j}{t_j - t'_j}(t_j - T_j) \geq Q'_j$$

controls the quality of activity "j" and this constraint can be written in the form:

$$\frac{q_j - q'_j}{t_j - t'_j}(t_j - T_j) \leq q_j - Q'_j$$

$Q'_j$ : Lower bound of quality for activity "j".

Assume we want to search for the optimal solution of cost, time, and quality at the same time such that:
$Q'_F = 0.98, Q'_I = 0.96$

Therefore, our model (17) will be in the form:

$$\begin{aligned}
& Max\ \lambda \\
& ST \\
& \lambda \leq \frac{1}{2} Th\left(\left(\frac{3060000 + 4250000}{2} - Z_1\right)\left(\frac{6}{4250000 - 3060000}\right)\right) + \frac{1}{2} \\
& \lambda \leq \frac{1}{2} Th\left(\left(\frac{42 + 29}{2} - Z_2\right)\left(\frac{6}{42 - 29}\right)\right) + \frac{1}{2} \\
& \lambda \leq \frac{1}{2} Th\left(\left(\frac{0 + 1}{2} - Z_3\right)\left(\frac{6}{1 - 0}\right)\right) + \frac{1}{2} \\
& t'_j \leq T_j \leq t_j \quad ; j = 1,\ldots,9 \\
& Y_i + T_i \leq Y_j \quad j = 1,\ldots,9 \\
& \lambda \geq 0 \\
& \frac{q_j - q'_j}{t_j - t'_j}(t_j - T_j) \leq q_j - Q'_j ; j = F, I
\end{aligned} \quad (18)$$

**The optimal solution is:**
$\lambda = 0.6133791$ ; $Z_1 = 3390000$, $Z_2 = 35$, $Z_3 = 0.41$
at
$X = \{T_A = 6, T_B = 8, T_C = 4, T_D = 4, T_E = 3, T_F = 7, T_G = 5, T_H = 6, T_I = 10\}$.

Note that, the optimal solution for model (17) was reduced into $\lambda = 0.6133791$ and the optimal results for the three criteria $Z_1, Z_2, Z_3$ changed as follows:

$Z_1: 3440000 \rightarrow 3390000, Z_2: 34 \rightarrow 35, Z_3: 0.34 \rightarrow 0.41$

It is expected that it will be changed because we will add constraints, but the changes are balanced for all criteria since as mentioned previously, the model was built to conduct many tests to get the best results for all criteria at the same time.

In another case, the decision maker can control the time of a project by adding constraint related to it. Suppose we want to find the optimal solution of cost, time, and quality at the same time. If we want to finish the project at most within 38 weeks, therefore we can add the following constraint $Y_{End} \leq 38$. There are many other cases in which our proposed model can absolutely help decision makers be in control and get the best solutions for the main project criteria: cost, time, and quality.

## V. DISCUSSION

The obtained results show the capability of the proposed model to obtain the optimal solutions for time, cost, and quality. Moreover, through using the membership function, the model is also able to control the quality factors for some important activities on which the project depends. This could be advantageous for decision makers in their initial planning stage. In addition to the fact that the results of the previous researches were similarly produced by our model, extensive experimentation showed that the proposed model is tolerant to add additional constraints related to each of the three criteria: time, cost, and quality.

## VI. CONCLUSION AND FUTURE WORK

Cost, time, and quality criteria for assessing projects form the solid triangle of project management. Previous research work used to solve this problem assigned a separate mathematical model for each criterion to find the optimal trade-off between the three of them. However in this paper, we used only one fuzzy mathematical model to solve this problem using non-linear membership function. Given the results of previous studies, numerical examples show that our proposed model is capable of solving the problem understudy and obtain the optimal solution of time, cost, and quality. However, further research and experimentation is useful for studying the effect of reducing a number of objective functions on increasing the tolerance of the proposed model for additional constraints. This could greatly positively impact decision making. Moreover, the study can be modified for further testing by adding other constraints according to the requirements of the problem, applying other combinations of the three criteria, or using other types of membership functions.

## REFERENCES

[1] Babu, A.J.G and Suresh. N., project management with time, cost, and quality–European J. Operations Research,-1996,320-327.

[2] Do Ba KHang, Yin Mon Myint, "Time cost and Quality trade-off in project management: case study", International Journal of Project Management, vol,17,No,1999,pp,249-256.

[3] Hamdy A. Taha.,2007-Operations Research: An introduction., University of Arkansas, Fayetteville,813.

[4] Kaheled El-Rayes, Amr Kandil, Time-cost-quality trade-off analysis for highway construction, Journal of Construction. Engineering Management, 2005 477-485.

Int. J. Sci. Res. in Computer Science and Engineering    Vol.**8**, Issue.**6**, Feb **2020**[5] Shrivastava. R, Singh.S and G. C. Dubey, India: Multi Objective Optimization of Time Cost Quality Quantity Using Multi Colony Ant Algorithm, Int. J. Contemp. Math. Sciences, Vol. 7, 2012, no. 16, 773 – 784.
[6] A. Afshar, A. Kaveh and O.R. Shoghli, Multi Objective Optimization of Time Cost Quality Quantity Using Multi Colony Ant Algorithm, Iran, Asian, journal of civil engineering vol. 8, no. 2 (2007)- 113-124.
[7] Kagade., K.L. and Bajaj., V.H,fuzzy method for solving multi-objective assignment problem with interval cost, journal of statistics and Mathematics,vol.1,issue1,2010,pp:1-9.
[8] P.K.De and DhartiYadav, An algorithm to solve multi –Objective assignment problem using interactive fuzzy Goal programming approach,Int.J.Contemp. Math.Sciences,vol.6,2011, no.34,1651-1662.
[9] P.Biswas And S. Pramanik, Multi-Objective Assignment problem with Fuzzy Costs for the Case of Military Affairs, International Journal of Computer Applications, vol 30-No,10,September 2011.
[10] BodkheS.G.,BajajV.H,Fuzzy Programming technique to solve bi-objective transportation problem, international journal of machine-intelligence, vol 2,issue 1,2010,pp:46-52.
[11] N. Ganapathy Ramasam, I. Mohamed Shaik Hameed, P. Ajai and Priya Sera Varkey,-November 2016, A Correlated Study between Time and Cost in Accordance with Fuzzy Logic, Indian Journal of Science and Technology, Vol 9(44), DOI: 10.17485/ijst/2016/v9i44/98800
[12] S. M. Mirmohseni, S. H. Nasseri, and A. Zabihi., 2017- An Interactive Possibilistic Programming for Fuzzy Multi Objective Solid Transportatio Problem, Applied Mathematical Sciences, Vol. 11, 2017, no. 45, 2209 – 2217
[13] Ehsan Ehsani, Nima Kazemi, Ezutah Udoncy Olugu, Eric H. Grosse, Kurt Schwindl., 2016- Applying fuzzy multi-objective linear programming to a project management decision with non-linear fuzzy membership functions, Neural Computing and Applications. . 10.1007/s00521-015-2160-0.
[14] J.Singh, fuzzy logics and fuzzy sets theory based multi-objective-Decision Making, Modeling Tool generalized for planning for Sustainable Water Supply, Institute of Town Planners, India Journal, 3,July-September 2011,73-80.
[15] kagadeK.L .and BajajV. H, Fuzzy approach with linear and some non-Linear membership functions for solving multi-objective assignment problems, advances in computational research, vol1,issue,2,2009 pp:14-17.
[16] P.K.De,A General Approach for Solving assignment problems Involving with Fuzzy Costcoefficients,vol,6,No,3:March 2012.
[17] S.Behera and J. Nayak, Solution of multi-objective Mathematical programming problems in fuzzy approach, vol 3 No,12 December 2011.
[18] Control: Managing Engineering, Construction and Manufacturing Projects to PMI, APM and BSI Standards,. Elsevier Science & Technology Books,431.
[19] Frederick S., Mark S., 2011- PERT/CPM Models for Project Management. In: Introduction to Management Science, Universityof Washington, PP:1-22.
_______________________________________________

### Appendix

```
max= k;
k>=0;
K<=0.5*@TANH((((3060000+3330000)/2)-
z1)*(6/(3330000-3060000)))+0.5;
K<=0.5*@TANH((((42+36)/2)-z2)*(6/(42-
36)))+0.5;
K<=0.5*@TANH((((0.3275)/2)-
z3)*(6/(0.3275)))+0.5;
tA<=10;tB<=9;tC<=4;tD<=6;tE<=3;tF<=7;tG
<=5;tH<=6;tI<=10;
tA>=6;tB>=7;tC>=3;tD>=4;tE>=2;tF>=5;tG>
=3;tH>=3;tI>=7;
yA>=0;yB>=yA+tA;yC>=yB+tB;yD>=yB+tB;yE>
=yB+tB;yF>=yC+tC;yF>=yD+tD;yG>=yD+tD;yG
>=yE+tE;yH>=yE+tE;
yI>=yF+tF;yI>=yG+tG;yI>=yH+tH;yEnd>=yI+
tI@;GIN)tA@ ;(GIN)tB@ ;(GIN)tC (
@;GIN)tD@ ;(GIN)tE@; (GIN)tF;(
-7)*0.035tF)<=1-0.9;8
-10)*0.066667tI)<=1-0.9;6
-10 )*0.0375)*(9/1)-1tA)+0.06*(9-
tB)+0.05*(4-tC)+0.1*(6-tD)+0.18*(3-
tE)+0.035*(7-tF)+0.05*(5-tG+(
-6)*0.083333tH)+0.066667*(10-
tI))>=0.96;
Z1=3060000+50000*( 10-tA)+50000*(9-
tB)+60000*(4-tC)+40000*(6-tD)+10000*(3-
tE)+40000*(7-tF)+75000*(5-tG(
-6)*90000+tH)+80000*(10-tI;(
Z2=yEND;
z3=0.0375*( 10-tA)+0.06*(9-tB)+0.05*(4-
tC)+0.1*(6-tD)+0.18*(3-tE)+0.035*(7-
tF)+0.05*(5-tG+(
-6)*0.083333tH)+0.066667*(10-tI;(
Q=1-(1/9)*(0.0375*( 10-tA)+0.06*(9-
tB)+0.05*(4-tC)+0.1*(6-tD)+0.18*(3-
tE)+0.035*(7-tF)+0.05*(5-tG+(
-6)*0.083333tH)+0.066667*(10-tI;((
```

@GIN)tG@; (GIN)tH@ ;(GIN)tI@; (GIN)yEnd (
@;GIN)yA@ ;(GIN)yB@ ;(GIN)yC (
@;GIN)yD@ ;(GIN)yE (
@;GIN)yF@;(GIN)yG@ ;(GIN)yH;(

@GIN)yI; (

```
end
```

**6**

**Authors Profile**

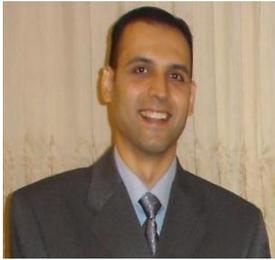

**Mohammad Mohammad Riad Sammany**
Faculty of Science – Department of Mathematics – University of Aleppo – Aleppo, Syria.

**Ph.D.** Applied Mathematics and Programming, Department of Mathematics, Faculty of Science, University of Aleppo, Aleppo, Syria, 2010
**Thesis Title:** "*Using Artificial Intelligence Techniques for Solving Fredholm Integral Equations and their Mathematical Applications*"
**Degree:** Distinction with Honors with an average of 98.4% (*the highest grade obtained by a post-graduate student in the Syrian Arab Republic Universities*)
**M.Sc.** Computational Mathematics, Faculty of Science, Department of Mathematics, Cairo University, Cairo, Egypt, 2007.
**Research Interests:**
•Using Artificial Intelligence Techniques for solving mathematical problems.
•Developing Numerical Methods for Solving Integral Equations and Ill-posed Problems.
•Developing Artificial Intelligence Techniques (Artificial Neural Network, Genetic Algorithms, Support Vector Machine, and Rough Sets) for solving real world applications.
**Honors and Awards:**
•2003: M.Sc. Governmental Fellowship for outstanding students, Ministry of Higher Education and Scientific Research, Syria
•2007: Testimony from Faculty of Science, Cairo University, Egypt, states that: *"The M.Sc. thesis of Mr. Sammany is quite sufficient to obtain a Doctorate Degree in Philosophy of Sciences"*.

• 2017: Full Post Doctorate scholarship granted by E.U. ERASMUS MUNDUS PARTNERSHIP PROJECT under the title "*Using Mathematical Models for Energy Saving in Buildings*". ACTION 2 Lot 1 – Syria *"IntegrAted Studies for Syrian and eUropean univeRsities"- ASSUR* Application reference number 551742-EM-1-2014-1-IT-ERA MUNDUS-EMA21

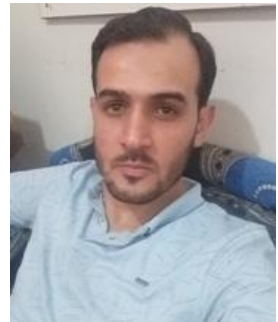

**Ahmad Steef**

**Ph.D.** Applied Mathematics – Informatics and Programming, Department of Mathematics, Faculty of Science, University of AL-Baath, Homs, Syria, 2017.
**MSc:** Applied Mathematics and Programming, Department of Mathematics, Faculty of Science, University of Aleppo, Aleppo, Syria, 2012.
**BSc.** Mathematics, Faculty of Science, Department of Mathematics, Tishreen University, Latakia Syria, 2009.
.

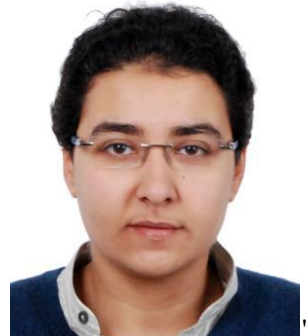

Nedaa Agami was a lecturer of Operations Research and Management Science (OR/MS) in the Department of Operations Research and Decision Support, Faculty of Computers and Artificial Intelligence at Cairo University. She is currently heading the Data Science and Advanced Analytics team in the Commercial International Bank (CIB). Prior to that, Nedaa worked as a Customer Intelligence and Data Mining consultant at Vodafone Egypt for seven years and has generally been offering consultancy and training services in the domain of Advanced Analytics since 2005. Nedaa received her PhD in January 2013 and her M.Sc. degree in 2008 in the same specialization. She has more than 14 years of experience in the field of Data Science, Advanced Analytics and Modelling in various industries including (but not limited to) telecommunications, tourism, health care, education, banking, supply chain, pharmaceuticals and many others for which she received several awards and has global best practices .

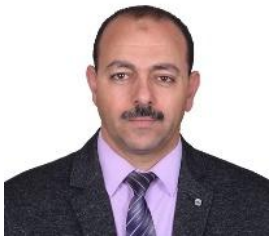

Tamer Medhat, born in Kafrelsheikh, Egypt, in July 13th, 1974, and received his Ph.D. degree in Faculty of Engineering, Tanta University, Tanta, Egypt in 2007. He is currently an Assistant Professor at the Faculty of Engineering and Vice-Dean of the Faculty of Artificial Intelligence for Education and Student Affairs, Kafrelsheikh University, Egypt. His current research interests include Information Systems, Artificial Intelligence, Augmented Reality, Decision Making, Computer Science, and Rough Set Theory Applications.